\documentclass{article}

\usepackage{amssymb}
\newtheorem{thm}{Theorem}

\date{}

\title{On the $p$-defect of character degrees of finite groups of Lie 
type}

\author{Meinolf Geck}

\begin{document}
\maketitle

\begin{abstract} \footnotesize 
This paper is concerned with the representation theory of 
finite groups. According to Robinson, the truth of certain variants of 
Alperin's weight conjecture on the $p$-blocks of a finite group would 
imply some arithmetical conditions on the degrees of the irreducible 
(complex) characters of that group. The purpose of this note 
is to prove directly that one of these arithmetical conditions is true in 
the case where we consider a finite group of Lie type in good characteristic.
\end{abstract}

According to Robinson \cite[Theorem~5.1]{Rob} (see also \cite[\S 5]{Rob2}), 
the truth of certain variants of Alperin's weight conjecture on the $p$-blocks 
of a finite group would imply some arithmetical conditions on the $p$-parts 
of the degrees of the irreducible (complex) characters of that group. The 
purpose of this note is to prove directly that one of these arithmetical 
conditions is true in the case where we consider a finite group of Lie type
in good characteristic. (See Example~2 for the problems which arise in bad 
characteristic.)

Let $G$ be a connected reductive group defined over the finite field
${\mathbb F}_q$, where $q$ is a power of some prime~$p$. Let $F \colon G
\rightarrow G$ be the corresponding Frobenius map and $G^F$ the
finite group of fixed points. Recall that $p$ is good for $G$ if $p$ is
good for each simple factor involved in~$G$, and that the conditions for 
the various simple types are as follows.
\[\begin{array}{rl} A_n: & \mbox{no condition}, \\
B_n, C_n, D_n: & p \neq 2, \\
G_2, F_4, E_6, E_7: &  p \neq 2,3, \\
E_8: & p \neq 2,3,5.  \end{array}\]
For the basic properties of finite groups of Lie type, see \cite{Ca2}.
Now we can state:

\begin{thm} Assume that $p$ is a good prime for $G$. Let $\chi$ be an
irreducible character of $G^F$. Then there exists an $F$-stable parabolic
subgroup $P \subseteq G$ and an irreducible character  $\psi$ of
$U_P^F$ (where $U_P$ is the unipotent radical of $P$) such that 
$|U_P^F|/\psi(1)$ equals the $p$-part of $|G^F|/\chi(1)$.
\end{thm}

In order to prove this result, we first reduce to the case that $G$ has
a connected center. This can be done using regular embeddings (see 
\cite{LuDis}), as follows.  We can embed $G$ as a closed subgroup into some 
connected reductive group $G'$ with a connected center, such that $G'$ has 
an ${\mathbb F}_q$-rational structure compatible with that on $G$ and $G,G'$ 
have the same derived subgroup. Then we have in fact ${G'}^F=G^F{T'}^F$ for 
any $F$-stable maximal torus $T' \subseteq G'$. In particular, the quotient 
${G'}^F/G^F$ is an abelian $p'$-group. Now, if we take any irreducible 
character $\chi'$ of ${G'}^F$ then, by Clifford's Theorem, the restriction 
of $\chi'$ is of the form $e(\chi_1 +\cdots +\chi_r)$ where $e\geq 1$ and 
$\chi_1,\ldots, \chi_r$ are irreducible characters of $G^F$ which are 
conjugate under the action of ${G'}^F$; moreover, $e$ divides $[{G'}^F:G^F]$. 
(In fact, Lusztig \cite{LuDis} has shown that $e=1$ but we do not need this 
highly non-trivial fact here.) It follows that both $e$ and $r$ are prime 
to $p$, and so the terms $|{G'}^F|/\chi'(1)$ and $|G^F|/\chi_i(1)$ (where
$1 \leq i \leq r$) have the same $p$-part. On the other hand, the collection
of unipotent radicals of $F$-stable parabolic subgroups is the same in $G$ 
and in $G'$. Thus, if the above theorem holds for $G'$, then it also holds 
for $G$.

Now let us assume that the center of $G$ is connected, and let $\chi$ be an 
irreducible character of $G^F$. We claim that there exists an $F$-stable 
unipotent class $C$ of $G$ such that \[q^{\dim {\mathcal B}_u} \quad 
\mbox{is the $p$-part of $\chi(1)$},\] where $u \in C$ and ${\mathcal B}_u$ 
denotes the variety of Borel subgroups containing~$u$. This class $C$ can be 
characterized in two different ways. One could either use the results in
\cite[Chap.~13]{LuBook} which describe a map $\xi_G$ from irreducible 
characters to unipotent classes in terms of the Springer correspondence. 
Then $C=\xi_G(\pm D_G(\chi))$ where $D_G$ denotes Alvis--Curtis--Kawanaka 
duality, and it remains to use the formula for $\chi(1)$ in 
\cite[(4.26.3)]{LuBook}. Or one uses the fact that every irreducible 
character has a unipotent support (see \cite{LuUS} for large~$p$, and 
\cite{Ge2} for the extension to good $p$). Note that here the assumption 
that $p$ is good is used, since then the denominator in the formula for 
$\chi(1)$ in \cite[(4.26.3)]{LuBook} is not divisible by $p$, which is 
definitely not the case in general. 

Using the theory of weighted Dynkin diagrams, one can associate with an
$F$-stable unipotent class $C$ two subgroups $U_2, U_1 \subseteq G$ such 
that 
\begin{itemize}
\item[(a)] $U_1$ is the unipotent radical of some $F$-stable  parabolic 
subgroup $P \subseteq G$, 
\item[(b)] $U_2$ is a closed $F$-stable subgroup of $U_1$ which is 
normal in $P$, 
\item[(c)] There exists some $u \in C \cap U_2$ such that the $P$-orbit
of~$u$ is dense in $U_2$ and $C_G(u) \subseteq P$.
\end{itemize}
If $p$ is large, this can be found in \cite[Chap.~5]{Ca2}, for example. It 
has been checked by Kawanaka \cite[Theorem~2.1.1]{Kaw2} that these results 
remain valid whenever $p$ is a grood prime. Now let $u \in C \cap U_2$ be as
in~(c). Then, by \cite[Lemma~2.1.4]{Kaw2}, we have 
\[ \dim C_G(u)=\dim C_P(u)=\dim P -\dim U_2.\]
On the other hand, we always have
\[ \dim P +\dim U_1=\dim G=2N+\mbox{rank}\, G,\]
where $N$ denotes the number of positive roots in~$G$. Using also the
formula $\dim C_G(u)=\mbox{rank}\, G+2\dim {\mathcal B}_u$ (see 
\cite[Cor.~6.5]{Spa}), we find that
\[ 2(N-\dim {\mathcal B}_u)=\dim U_1+\dim U_2.\]
Now Kawanaka \cite[\S 3.1]{Kaw2} has shown that there exists a linear 
character $\lambda$ of $U_2^F$ which is invariant under $U_1^F$ and such that
$U_1^F/\mbox{ker}\lambda$ is an extraspecial $p$-group. By the character
theory of extraspecial $p$-groups, $\dim U_1 -\dim U_2$ is even and there 
exists an irreducible character $\psi$ of $U_1^F$ which lies above 
$\lambda$ and has degree $q^{(\dim U_1-\dim U_2)/2}$. (Here, we use
that if $U$ is any connected unipotent group defined over ${\mathbb F}_q$,
then the fixed point set under the corresponding Frobenius map has
order $q^{\dim U}$.) Using the above formulas this yields  that 
\[ \frac{|U_1^F|}{\psi(1)}=q^{\dim U_1-(\dim U_1-\dim U_2)/2}=
q^{(\dim U_1+\dim U_2)/2}=q^{N-\dim {\mathcal B}_u}.\]
Since $q^N$ is exactly the $p$-part in $|G^F|$, the proof of Theorem~1 is
complete.

\medskip
\noindent {\bf Example~1}. Let $G^F=\mbox{GL}_3(q)$ and $\chi$ be of degree 
$q(q+1)$. The $p$-part of $|G^F|/\chi(1)$ is~$q^2$. The unipotent class $C$ 
associated with $\chi$ has Jordan blocks of sizes $1,2$. The group $U_1$ is 
the group of all upper triangular unipotent matrices, $U_2$ is the center 
of~$U_1$, and for each nontrivial linear character of $U_2^F$ there is a 
unique irreducible character of $U_1^F$ lying above it. Any such irreducible 
character of $U_1^F$ has degree~$q$.

\medskip
\noindent {\bf Example~2}. Let $G^F=\mbox{Sp}_4(q)$ and $\chi$ be of degree 
$\frac{1}{2}q(q^2+1)$. Then the unipotent support of $\chi$ is the unique 
unipotent class $C$ with $\dim {\mathcal B}_u=1$. (This holds independently 
of whether $q$ is a power of a good prime or not; see \cite{GeMa2}.) In good 
characteristic, the subgroups $U_1,U_2$ associated with $C$ are such that 
$U_1=U_2$ and $|U_1^F|=q^3$. Analogous subgroups will also exist in 
characteristic~$2$, but we cannot work with them since then the $2$-defect 
of $\chi$ is $2q^3$. --- This illustrates some of the difficulty of 
extending Theorem~1 to the case of bad characteristic.

\medskip
\noindent {\bf Acknowledgements.} This note was written while I was 
participating in the special semester on representations of algebraic
groups and related finite groups at the Isaac Newton Institute (Cambridge,
U.K.) from January to July 1997. It is a pleasure to thank the organisers,
M.~Brou\'{e}, R.W.~Carter and J.~Saxl, for this invitation and the Isaac
Newton Institute for its hospitality. 

%%%%%%%%%%%%%%%%%%%%%%%%%%%%%%%%%%%%%%%%%%%%%%%%%%%%%%%%%%%%%%%%%%%%%%%%%%%

{\sc\footnotesize Institut Girard Desargues, Universit\'e Lyon 1,
21 av Claude Bernard, F--69622 Villeurbanne cedex, France}

\medskip
\makeatletter
{\footnotesize\it E-mail address: \tt geck@desargues.univ-lyon1.fr}
\makeatother
\end{document}